\documentclass[10pt,twoside]{siamltex}
\usepackage{amsfonts,epsfig}

\usepackage{xy}
\input{xy}
\xyoption{all} \xyoption{rotate}

\newtheorem{example}[theorem]{Example}

\begin{document}

\bibliographystyle{plain}
\title{
On coordinatization of mathematics }

\author{
Peteris\ Daugulis\thanks{Department of Mathematics, Daugavpils
University, Daugavpils, LV-5400, Latvia (peteris.daugulis@du.lv).
} }

\pagestyle{myheadings} \markboth{P.\ Daugulis}{On coordinatization
of mathematics} \maketitle

\begin{abstract} The problem of advancing logic based coordinatization of mathematics is
considered. The need to develop a theory for measuring value and
complexity of mathematical implications and proofs is discussed
including motivations, benefits and implementation problems.
Examples of mathematical considerations for such a theory are
given. Arguments supporting applications in mathematical research
guidance, publication standards and education are given.

\end{abstract}


\section{Introduction}\

\subsection{The aim of the article} The main aim of this article is to
point out the need and possibilities to interpret and encode
logical implications (inferences, consequences) and proofs as
mathematical objects focusing on complexity and paying less
attention to semantic specifics and syntactic issues. Advancement
of coordinatization of mathematics is meant in the sense of
finding the simplest mathematical structure faithfully
representing mathematical proofs, theories, creativity and
development. Motivations and possible benefits of this idea are
discussed in a programmatic style. This idea must lead to
important advances such as proof complexity measures and models
for networks of mathematical statements. It must have important
applications in research guidance, evaluation of mathematical
results and education, these applications seem to be important in
their own right. Arguments which show that the proposed program
will have new features compared to classical logic, computational
logic (automated theorem proving), category theory and
computerized projects related to formalization of mathematics are
mentioned. There are no theorems in this article, most issues are
discussed with a certain vagueness. This article is not intended
to contribute to literature related to an established problem
although one can find features reminding of the $24$th Hilbert's
problem. It can be interpreted as both a research proposal and a
discussion oriented report. The paper is oriented mainly towards
mathematicians interested in logic and philosophy (expansion of
knowledge) of mathematics.

\subsection{History and the current state}
The usefulness and the exceptional role of
mathematics have been best expressed by the hypothetical
Pythagorean saying "all is number" which asserts that all physical
objects, systems and processes may be precisely mathematically
modelled using numerical constructions. Using a philosophical
point of view and language mathematics can be thought as a
universal epistemological framework created by the human intellect
in order to justify knowledge or to perform justification/regress
steps, see Pollock (1975), for knowledge from various areas in a
uniform way.

%

Mathematical activities have produced both pure and applied
results: ma\-the\-ma\-ti\-cal theories - collections of related
mathematical facts, algorithms and proof techniques, nontrivial
properties of mathematical objects, nonobvious logical
con\-seq\-uen\-ces and practical computational algorithms.

The progress of both pure and applied mathematics has been greatly
influenced by advances in formalization and coordinatization of
mathematical objects and the mathematical language. Encoding
techniques of mathematical primitives, objects, statements and
computational steps have influenced the progress of all sciences
and most human activities. As encoding breakthroughs in
mathematics one can mention the introduction of the positional
notation, al\-geb\-ra\-ic operations, Cartesian coordinates,
calculus, matrix algebra, mathematical logic, category theory etc.
See Shafarevich (1990) for a treatment and examples of
coordinatizations in algebra.  Many syntactic problems of
mathematical statements have been considered and solved.
Significant areas of mathematical logic such as the first-order
logic, Gentzen and Hilbert calculi etc. were founded, developed
and applied. For some areas of mathematics minimal systems of
axioms and inference rules have been proposed, for example, the
Peano arithmetic. The term ``metamathematics'' was introduced by
Hilbert to denote advanced mathematical logic. Recently there has
been an attempt to do mathematical logic without syntax - to
define and study \sl combinatorial proofs\rm\ in propositional
logic as graph homomorphisms of certain kind, to distinguish
between \sl syntactic witnesses\rm\ and \sl mathematical
witnesses\rm\ of proofs, see Hughes (2006).

Another important breakthrough for the mathematical language and
thinking was the category theory introduced by Eilenberg and Mac
Lane. The category theory is a successful attempt to advance
unified encoding of the mathematical language and establish maps
between different mathematical theories.

One can notice two major longterm trends in mathematics and its
applications. First of all, the application areas having precise
mathematical models and being served by applied mathematics are
constantly enlarging. Complexity of mathematical models is
steadily increasing, for example, consider models of biological
systems and processes. Even the most ``unmathematical'' notions
and processes, for example, related to consciousness and
psychological activities, may be subject to mathematical modelling
(in particular, due to emergence of mathematical models of nervous
systems) - justification/regress steps in the philosophical sense.
Secondly, mathematical notations and the mathematical language are
steadily getting more rigorous and well defined over time. We call
these two trends \sl the Pythagorean process.\rm\


Computers are constantly being used by mathematicians also to
 check
and verify theorems, to make computations termed ``automated
theorem proving`` and "mechanical theorem proving" which are
equivalent to statement proving in certain areas of mathematics,
see Robinson (2001) and Chou (1994). Proposals for computer based
mathematical knowledge management systems such as ``QED
manifesto'', see Wiedijk (2007), have been made. See Avigad,
Harrison (2014) for a recent review of computerized theorem
checking/proving.

\subsection{Possible nexts steps}

Advances of the Pythagorean process have always been significant
events in human history shaping sciences, technologies and
thinking. One can ask whether the Pythagorean process will
continue and what may be its next steps. In order to predict and
stimulate nontrivial and important advances we must look at those
features of mathematical activities which have not been
coordinatized and measured yet. The goal of this paper is to show
and discuss a possible direction for the Pythagorean process.


One feature that is still missing in the mathematical culture is
modelling and representation of creative mathematical thinking
going beyond its semantic and syntactic content - a precise
expression of all mathematical implications, proofs and algorithms
as well defined mathematical objects. Going beyond semantic and
syntactic content means development of implicational
propopositional calculus - mapping these features to universal,
simpler mathematical objects. Another missing feature is
mathematical structure and representations of mathematical
theories - collections of related mathematical results and proof
techniques. In philosophical terms these missing features
correspond to a justification/regress step with respect to the
mathematics itself - a justification of mathematical implications
in tems of relatively simple mathematical concepts.

We conjecture that there will be new developments of mathematics
which will continue the Pythagorean process - there will be
advances of the mathematical thinking and encoding which will
allow us to go beyond semantics and syntax of mathematical texts -
to comprehensively coordinatize (map into mathematical objects)
and precisely interpret proofs and ma\-the\-ma\-ti\-cal theories
as numerical or geometrical mathematical objects, known or new
ones. Given a mathematical theory $A$ (a structure containing
objects of study, first-order or higher-order logic statements,
proofs etc.) we may look for a mathematical object $\alpha$ which
would be a good model of $A$: elements of $A$ such as logical
implications, proofs and subsets of mathematical statements in $A$
would be defined as substructures or quotient structures of
$\alpha$. The transfer from $A$ to $\alpha$ should be thought
philosophically as a regress step.

The proposed idea goes beyond the standard mathematical logic
which deals with constructions of systems of axioms, correct
statements, syntactic and language problems, expressive power
problems of axiom and inference systems. The proposed research
program also goes beyond programs such as Hilbert's program and
the recent ``QED manifesto'' program because of its focus on
implications and models of theories. Our idea can be roughly
compared to introducing Cartesian coordinates - assigning
implications directions and lengths whereas the standard classical
logic is interested mainly in premises and conclusions. For the
same reason it goes beyond the computational mathematical logic
(e.g. automated theorem proving) which deals with computerized
proving or disproving statements in a given formal language.
Programs for creation of computerized data systems of mathematical
knowledge have not been completely successful because they did not
focus on implications. Our proposed program may have links to
proof complexity theory. Our idea may also provide new examples
and cases of study for the category theory.

Such models would allow to increase the speed and improve the
quality of progress of a given mathematical theory, improve
understanding of various theories, compare different theories and
improve understanding of their relations, measure and quantify
mathematical results such as theorems and lemmas, classify
mathematical theories up to isomorphism in a right sense, consider
and interpret maps between mathematical theories. It would enable
mathematicians counteract the specialization drive and digest
bigger amounts of information. This research proposal is related
to the lesser-known $24$th Hilbert's problem - find the simplest
proof of a given statement, compare different proofs, design
criterions for simplicity and rigor etc., see Thiele (2003).
Finding mathematical models of proofs should be considered the
main unsolved problem in mathematics nowadays containing the
$24$th Hilbert's problem as a subproblem.

These models may also provide new research questions about the
object $\alpha$ or stimulate development of new mathematical
objects for the described modelling purposes. Different models may
be initiated by different areas of mathematics, these different
models may formalize and extract specific research experiences.
They may provide one more abstraction step in human thinking -
allow to make logical implications without focusing on semantic
content of premises and conclusions. From the computational point
of view it may allow to substitute logical implication making by
computations. The coordinatization should be essentially unique.
It would be considered as a technique for reducing complexity and
creativity needed for work in an area, i.e. the mathematical
research, to a minimal admissible level. It can also lead to
generalizations of the implication concept.

Such a development would allow to describe the current state and
development, in narrow or wide areas, and maybe even the whole
history of mathematics as a mathematical object, to measure
mathematical proofs, algorithms and results in a mathematically
well defined way, to classify and measure mathematical creativity,
to derive canonical solution paths for unsolved problems. It would
also be used to guide researchers, show them the most important
research directions, problems and milestones in a rigorous and
quantitative way. Both theory building and problem posing/solving
have to be formalized. The mathematical creativity, the progress
of mathematics itself and the goal of mathematics has to be
defined as mathematical objects. Such an advance of the
Pythagorean process may generate new encodings and metalanguages
for mathematical statements and proofs. Its successes in pure
mathematics may be transferred to other sciences through applied
mathematics and thus mathematics one more time may play a decisive
role in human history.

If future generations will be interested in further mathematical
research (especially in pure mathematics) then computers or their
future descendants will be eventually used to perform it.
 Therefore we need to create theories which would interpret and
model human mathematical thinking using mathematical objects which
can be processed by computers, reduce mathematical goal setting
and creative theorem proving to computation, define the goal of
mathematics as a computational result. The step of passing from
computations to proofs and algorithms should be iterated producing
new paradigms of transformations of proofs and algorithms. It may
be impossible to change human thinking but it may be realistic to
organize and emulate a mathematical research process which would
be performed by computers or future virtual minds. Even if this
project is not successful mathematical thinking has to be
evolutionized in order to take into account computing
technologies.

Results of implication coordinatizations and modellings will
advance our understanding of implication making and thinking
itself to new levels, question the role and the very need of
implication making, offer possible improvements. It may carry
mental and computational activities to new heights, identify
limitations, weaknesses and pecularities of human thinking. If
this approach is successful we may ask fundamental questions: what
can be considered an advanced or future form of mathematical or
general implication/consequence making? if there is such a form
how it can be implemented? when will human made computing devices
outdo humans in creative mathematical implication making and what
will be the consequences?

Although it is not within the scope of this paper it can be
mentioned that possible results in the proposed direction may be
combined with expected advances in biology related to detailed
description and understanding of the organization of the brain
functioning on the subcellular level. Mathematical thinking as
cognitive activity, language processing and thinking in general
must be modelled, analyzed and modified, if possible, starting
from the detailed analysis of the natural physical process in the
brain.

There may already exist scattered examples which are known to
experts and the Pythagorean process may procceed in the proposed
direction spontaneously. Ne\-ver\-theless relevant results and
examples should be integrated into a single program. Even if the
proposed research projects are not considered successful and
results are negative the related work may have partial successes -
it may generate nontrivial mathematical results, higher levels of
abstraction, new encodings and standards for mathematical
language, discourse and objects, new thinking features and
philosophy. We have to allow the possibility that asking a
question may be more important than solving it.

\subsection{Applications}

A mathematically sound method for measuring value or complexity of
mathematical results would also allow to set rigorous standards
for research publications in profesionally and internationally
accepted journals and other information depositories. Currently
mathematical results are evaluated apparently without any rigorous
system. Additionally, in some areas applied mathematics is used
due to its "decorative" value. The current competition oriented,
trend based and partially dogmatic evaluation  of results and
merits can not be considered justified in mathematics which is the
very source and center of the culture of unbiased logical
reasoning and numerical analysis.
The lack of a rigorous evaluation theory is a sign of backwardness
in the same way as the lack of mathematical modelling is such a
sign in any other area. A rigorous evaluation method based on
mathematical analysis of results and techniques must be found.


In section \ref{1} we give descriptions of these and other
possible applications.

\section{Main research and application directions}\

\subsection{Coordinatization of implications and proofs}

%
%
%
%
%

Logical foundations and consistency of mathematics have been
intensively studied since the design of the Hilbert's program, see
Simpson (1988). Decidability of mathematical statements and
inherent limitations of axiomatic mathematical systems have been
investigated since G\"{o}del, see Davis (2006). The first
G\"{o}del incompleteness theorem is obtained using arithmetization
- encoding of Peano axioms using natural numbers,
see Kozen (1997, pp.206-292). Thus G\"{o}del results are examples
of successful instances of coordinatization and modelling of
mathematical theory and language: encoding of mathematical
statements, modelling the set of mathematical statements using set
theory and interpreting set-theoretic results in terms of
statements and proofs.

Mathematical structure of theories, complexity of implications and
mathematical proofs and precise value of mathematical results
regardless of their semantics and syntactic issues have not been
paid adequate amount of attention of wide mathematical community.
Automated theorem proving of some first-order and second-order
logic statements still does not substitute the human intellect.
Apparently nontrivial logic - mathematical implication making and
creativity are considered to be areas which can not be
mathematically modelled and computerized. Another reason for this
may be that the natural biological limit of human intelligence has
been reached - human brain can not sustainably operate at such a
complexity and abstraction level. A limiting state will be reached
simply because human brain has limited physical capabilities
(number of nerve cells, neuron structure and neuron connection
graph complexity) and more complex computational devices will be
designed. We must prepare for such a limiting state. It may be
noted that mapping of mathematical proofs to simpler mathematical
objects may not necessarily involve higher abstraction levels.

Proofs of mathematical statements are sequences or, more
generally, networks of logical implications. Therefore one
approach to the study of proofs would be to study relatively
simple logical implications and their networks. Research may also
be needed to determine right definitions of irreducible
implications, various types of implications and their linkings,
embeddings of the objects corresponding to implications in
suitable ambient spaces - a geometrization of logic, definitions
of creativity. The concept of proof may be generalized for
undecidable statements. For such statements an analogue of proof
may be an infinite process converging in a right sense and
infinitary logic may need to be used.

Logical implications can be defined as operations on logical
predicates in first-order or higher-order logic using logical
connectives, especially the material condition connective
$\Rightarrow$. The consequence relation $\vdash$ used in
mathematical logic is also a relevant notion. Given two predicates
$P(x)$ and $Q(x)$ defined for all $x\in X$ we say that \sl $P$
implies $Q$\rm\ ($P\rightarrow Q$) provided
$$\bigwedge_{x\in X}\Big(P(x)\Rightarrow Q(x)\Big)=true.$$ The support
$supp(A)$ of a predicate $A$ may be defined as the set of $A$
argument values $x$ for which $A(x)=true$, thus $supp(A)\subseteq
X$. Validity of a predicate implication $P\rightarrow Q$ is
equivalent to the set-theoretic inclusion of the support of $P(x)$
into the support of $Q(x)$: $P \rightarrow Q$ is a true statement
if and only if $supp(P)\subseteq supp(Q)$.  We could try to
coordinatise the implication $P \rightarrow Q$ by set-theoretical,
combinatorial, algebro-geometrical, geometrical, topological and
complexity-theoretical properties of the sets $supp(P)$ and
$supp(Q)$ such as 1) absolute and relative sizes and shapes of
$supp(P)$, $supp(Q)$ and $supp(Q)\backslash supp(P)$, 2)
properties of the boundaries of $supp(P)$ and $supp(Q)$. We
conjecture that 1) the implication $P \rightarrow Q$ can be
considered easy if $supp(P)$ is a relatively small, e.g.
low-dimensional, subset of $supp(Q)$; 2) implications $P
\rightarrow Q_{1}$ and $P \rightarrow Q_{2}$ can be considered
distinct in a proper sense if $(supp(Q_{1})\cap
supp(Q_{2}))\backslash supp(P)$ is relatively small.

Proofs as sequences of implications $P_{1}\rightarrow
P_{2}\rightarrow ... \rightarrow P_{n}$ may be considered as
sequences of set-theoretic inclusions $supp(P_{1})\subseteq
supp(P_{2})\subseteq ... \subseteq supp(P_{n})$. Passing from
semantic-specific implication making to constructing sequences of
embedded sets should be considered as a computational substitution
of implication making.

Coordinatization and measurement of logical implications may also
be related or even reduced to computational complexity if
computations are involved determining the inclusion
$supp(P)\subseteq supp(Q)$.

Viewing proofs as directed paths in a proof graph or other ambient
structure with edges corresponding to simple implications we can
try to use the graph-theoretical or topological intuition
describing properties of proofs.

Additional idea is to generalize implications, to define other
binary relations in statement sets or consider weighted graphs
(not to be confused with fuzzy logic). In terms of generalized
implications standard implications would be their special case.
Given two predicates $P(x)$ un $Q(x)$ we can consider another
properties of sets $supp(P)$ and $supp(Q)$ (instead of inclusion)
for this purpose. For example, we can define that $P$ \sl almost
implies\rm\ $Q$ if $supp(P)\backslash supp(Q)$ is relatively small
or simple in a suitable sense.

\begin{example} We give a candidate definition for irreducible
implications in the case of propositional logic. Suppose
$p(X_{1},...,X_{n})$ and $q(X_{1},...,X_{n})$ are formulae in
propositional Boolean variables $X_{1},...,X_{n}$ and the
implication $p\rightarrow q$ is true. We call the implication
$p(X_{1},...,X_{n})\rightarrow q(X_{1},...,X_{n})$ irreducible if
the full disjunctive normal form (DNF) of $q$ has exactly one more
disjunctive term than the full DNF of $p$. The implication
$p\rightarrow q$ is not a composition of two noninvertible
implications.
\end{example}

\subsection{Modelling approaches}

Hilbert's point of view of proofs as mathematical objects should
be developed further. In this section we consider a few possible
directions for advancing Hilbert's ``logical arithmetic''.

\subsubsection{A category-theoretic approach to modelling of implications
and proofs}

A general approach for modelling mathematical theories and proofs
is category theory. Define a category $\textbf{Math}$ where
objects are mathematical statements and morphisms are logical
implications (consequences), composition of morphisms may be the
standard composition of implications. Subcategories of
$\textbf{Math}$ would correspond to specific mathematical theories
and functors between these theories would show their mappings. We
can try to study $\textbf{Math}$ or its subcategories with respect
to problems such as concretizations, functors to and from other
categories, interpretations of category-theoretic constructions
such as natural transformations, adjoint functors, pushouts and
pullbacks, limits, quotients etc.

\subsubsection{Graph-theoretic modelling of mathematical
theories}\label{g}

As we no\-ted in a previous section a mathematical theory can be
interpreted as a directed graph corresponding to the implication
relation which we call \sl proof graph\rm\ $\Pi=(\Sigma, \Lambda)$
with vertices in the set $\Sigma$ being statements (which are not
interpreted as implications) and directed edges in the set
$\Lambda$ being relatively simple logical implications although
the complexity of these implications depend on future models and
is a matter of further study.  Graph theory, see Diestel (2010),
can be considered as a candidate theory whose concepts and methods
may be used to coordinatise, measure, compare and visualize proofs
and mathematical theories. Although graphs are already widely used
in mathematical logic further in this subsection we give a few
examples of graph-theoretical considerations which may be useful.

%

\paragraph{Metric properties of the proof graph} Assume that any
edge of a proof graph $\Pi$ is given a weight which measures the
complexity or some other well defined property of the
corresponding implication. In the simplest naive cases weights
could be positive numbers but other weight sets can not be
excluded from consideration. Assume that we are given a directed
path between two vertices $P$ and $Q$ having edges
$e_{1},e_{2},...,e_{n}$ with weights $w_{1},w_{2},...,w_{n}$ which
corresponds to a proof $P \rightarrow Q$. Complexity or other
measure of the proof could be defined as an appropriate function
of weights $w_{1},w_{2},...,w_{n}$, for example, the sum
$w_{1}+w_{2}+...+w_{n}$. Most likely, more complex weight
functions dictated by mathematical logic will be used to describe
and classify implications and proofs. For example, a weight
function may assign each edge the corresponding implication type
in an appropriate sense. Having a proof graph invariant which
would correspond to proof weight or metric we could investigate
problems such as, for example, the problem of finding all
statements within a fixed distance from a given statement or
axiom. Analogs of various metric-based subgraphs such as nearest
neighbour graphs can be studied.

\paragraph{Vertices with special/extremal properties as valuable or unvaluable statements}

Proof graph models and other proof coordinatization ideas should
rigorously identify extremal relations, operations, statements and
extremal implication steps which are relatively more or less
important than others.

In particular, vertices of proof graphs having extremal properties
related to connectivity, metric, centrality or other invariants
may be considered as valuable "theorems". For example, the notion
of graph center could be suitably modified and vertices having
minimal weighted eccentrity defined as valuable statements. The
same arguments should identify statements which can be considered
of low value.

\paragraph{Path systems} Different paths in the proof graph between vertices
$P$ and $Q$ represent different proofs between the corresponding
statements. Having fixed vertices $P$ and $Q$ we can study all
$(P,Q)$-paths, e.g. we can pose the problem of finding all
$(P,Q)$-proofs up to a certain equivalence relation. We can also
try to find vertices with special properties, e.g. vertices which
are in more than one $(P,Q)$-path. In topological models for proof
spaces topological ideas such as homotopy classes of path systems
and homology-type invariants can not be excluded from
consideration.

\paragraph{Shortest paths} Given two statements $P$ and $Q$ in a
proof graph we could look for $(P,Q)$-paths with some special or
extremal properties such as the paths having minimal weight. That
would correspond to finding $(P,Q)$-proofs with some special
properties, for example, the proof of minimal complexity. These
ideas again remind us of the $24$th Hilbert's problem and the
``simplest proof''.


\paragraph{Sources and sinks} Strong equivalence classes of vertices of the proof graph having
no incoming or outgoing edges can be interpreted as axioms and
terminal statements, respectively.

\paragraph{Graph homomorphisms} Given two proof graphs we can
consider maps between them which preserve desired graph properties
similarly to graph homomorphisms or isomorphisms. Such maps,
perhaps linked with category-theoretic constructions, would allow
to define maps between corresponding theories, to compare and
classify theories, construct mathematics as a single object.

%
%

\subsubsection{An algebraic approach to modelling of implications
and proofs}

A mathematical theory can also be interpreted as an algebraic
structure as follows. Given two directed adjacent implications
$f:P\rightarrow Q$ and $g:Q\rightarrow R$ their composition
$g\circ f:P\rightarrow R$ is an implication. The composition of
implications can be interpreted as a binary associative operation
on the set of implications. Additionally the operation has to be
defined for nonadjacent implications. The implication set
$\Lambda$ thus has a na\-tu\-ral monoid structure $(\Lambda,
\circ)$, algebraic questions may be asked and algebraic methods
may be used to study $\Lambda$. For example, submonoids, ideals,
congruence relations and quotient structures of $\Lambda$ could be
studied and interpreted. In this approach we also may consider
generalizations of the implication notion.

\subsubsection{A topological approach to modelling of implications
and proofs} A mathematical theory $(\Sigma, \Lambda)$ can also be
endowed a topological space structure as follows. We start with
noting that the implication binary relation $\rightarrow$  is a
preorder relation - it is obviously reflexive and transitive. We
can view the implication relation as a specialization preorder for
the Alexandrov topology $\tau$ on $\Sigma$ corresponding to
$\leftarrow$: the open sets for $\tau$ are the upper sets with
respect to the relation $\leftarrow$. We remind the reader that a
set $U$ is an upper set with respect to $\leftarrow$ provided
$Q\in U$ and $P\rightarrow Q$ implies $P\in U$, see Barmak (2011).
Thus we can investigate the given mathematical theory $(\Sigma,
\Lambda)$ using topological experience and intuition - study the
topology $\tau$ with respect to standard problems of general and
algebraic topology such as interpretations of continuity or
(co)homology invariants.

\subsubsection{Complexity-theoretic approaches} Given an
implication or a proof $f:P\rightarrow Q$ we can measure the
(deterministic) complexity of $f$ as some computational complexity
measure (time or space related) of a computational process
producing $f$ with given $P$. An example of such measure can be
proof size considered in proof complexity branch of proof theory.
Nondeterministic complexity - computational complexity of proof
checking should be considered as well.

Value of mathematical results can be estimated considering their
impact on computation complexities (time, space, parallelability
etc.). A result can be considered valuable if it has a
computational value such as reduction of complexity classes of
computational and decision problems. On the contrary, a result may
be considered easy if it amounts to a polynomial time reduction.
History of mathematics should be studied as a network of
complexity reductions. Mathematical questions about structure of
reductions as maps between languages should be studied.

\subsubsection{Global analysis of mathematical results and theories}
Non\-tri\-vi\-al math\-ema\-ti\-cal results can be analyzed as
single objects and initial proofs can be standartized, improved
and optimized. For example, results can be analyzed with respect
to existence of Noetherian induction proofs. Suppose the statement
$\forall\ x\in X\ P(x)$ is true, does there exist a relatively
simple well-founded relation $R\subseteq X\times X$ such that the
statement can be proved relatively easy using Noetherian
(structural) induction on $R$ ? Complexity of involved
well-founded sets and induction steps can be considered as global
complexity and value measures.

\subsubsection{Proof bundles} If we have two predicates $P(x)$, $Q(x)$ where $x\in X$
and an implication or proof $f:P\rightarrow Q$ which is true for
every $x\in X$ then the complexity of proofs and proofs themselves
may be different for different $x\in X$. Such situations may be
considered using topological analogy with topological bundles, the
set $X$ being the base and the proof $f_{x}$ for each $x\in X$
being the fiber.

\subsection{Applications - research guidance and requirements for mathematical
texts}\label{1}

Mathematical results in form of correct proofs are described in
research papers and other documents such as monographs. Research
problems and new mathematical objects are often insufficiently
motivated. Some mathematicians seem to value a result just because
it describes an object, a property or a case which has not been
described before, this resembles publishing a computation result
just because it has not been published before and should not be
accepted in mathematics. Many ma\-the\-ma\-ti\-cians seem to
prefer research-like activities on insignificant problems or cases
which would guarantee a publishable result instead of working on
important hard problems which may be hopeless. Since the research
merit of a mathematician is effectively determined as a weighted
sum of numbers representing her/his publications in indexed
journals there are mathematical publications in these journals
which seem to have been published only to increase the number of
publications of their authors. There is no rigorous method to
determine the value and the nontriviality/creativity of a
mathematical result as a mathematical object. An advance of the
Pythagorean process is needed.

\subsubsection{Research guidance}

Mathematical research processes, problems, conjectures and
research interests should be motivated by rigorous analysis based
on a proof and statement coordinatization theory. Such a theory
would show meaningful problems, computations and/or directions
which need to be studied in order to advance the understanding of
a given domain, missing or optimal concepts that need to be
introduced, proofs that need to be modified, mathematical regress
steps (mappings to simpler objects) that need to be done etc. It
would direct and enforce development of mathematics, link and rank
various areas of mathematics more closely and effectively than
category theory. The progress of mathematics must be defined
mathematically and research must be performed in an optimal way.
Finding and studying motivations for research problems should be
part of research. The global ma\-the\-ma\-ti\-cal research process
should be more clearly than now subdivided into a number of
subprocesses corresponding to well defined longterm problems. As
examples of such longterm problems one can mention longstanding
conjectures in number theory and various classification problems
in algebra and topology. Alternatively, terminal states and
ultimate goals of all significant domains of mathematics should be
envisioned, roadmaps and design of ultimate problems (such as
classification, decidability or complexity-theoretic problems) for
reaching these states should be part of research. Complexities
involved in problem formulations should be compared to
complexities of their solutions. Problems which are easy to
formulate but difficult to solve should be highly valued and
studied for this reason.

\subsubsection{Standardization of mathematical proofs and texts} A
standard proof format could be designed for mainstream
mathematical publications. For example, the description of a proof
would be organized so that the steps are uniform and clearly
shown.  Standart layout formats for monographs and textbooks
showing the structure of mathematical theories and complexities of
each proof step would improve quality and clarity of mathematical
texts.

\subsubsection{Historical research} Although active mathematicians are mainly interested in unsolved problems and theory
building it makes sense to pose problems related to the history of
mathematics. Existing mathematical texts could be exhaustively
studied with respect to proof complexity, development of concepts
and proof structure. Proof structure of mathematical results which
are considered important should be analyzed. More generally,
events which have formed mathematics should be identified
rigorously even though they are already known. An ambitious
longterm collective goal in this direction for researchers could
be an exegesis and mapping of the body of most known and
reasonably important mathematical results produced by the mankind
to a suitable mathematical object which would record objects of
study, statements, proofs and the process of their production.

\subsection{Control of the publishing process} Apart from guiding
mathematical research and improving mathematical texts new
advances in proof coordinatization and complexity theory could
control the flow of published mathematical texts and formation of
research merits of mathematicians. Moreover, mathematical
modelling and analysis of the publication process, reforms of
mathematical discourse and professional development (both
individual and collective) should be an important part of
professional duties of logicians and philosophers.

\subsubsection{Reviewing}

Currently there is a large number of scientific journals, internet
based archives and other publishing opportunities. In most
journals the value and the originality of a correct new
mathematical result submitted to a journal is vaguely determined
as an emotional (not being rigorously determined) opinion of one
or more reviewers or even just the responsible editor of the
journal. No rigorous method to determine values of research
results seems to be known and used. Since the amount of published
mathematical results is steadily increasing it may be difficult
for a reviewer to determine the value of a new result. These
values are often determined by popular opinions in wide or narrow
scientific communities at the given moment of time. The appearance
of documents like \sl Code of Practice\rm\ created inside the
European Mathematical Society, see European Mathematical Society
Ethics Committee (2012), indicates a necessity of changes in the
research result value determination process. Mathematical texts
submitted for publication in journals or internet archives should
be analyzed much more rigorously and openly than it is done now.

Currently most mathematicians seem to believe that research
process can not be mathematically modelled, computerized and,
therefore, rigorously measured. In the author's opinion, this view
is wrong, unreasonably epistemologically conservative and
exceedingly focusing on the previous experience. This view of the
research process resembles performing arithmetical operations
before introduction of a numeral system.

A rigorous proof complexity and value theory would allow to define
and determine values of correct submitted or published results
more rigorously and set standards for them. Values of research
results may need to be defined locally (considering the state of
mathematics in a relatively short period of time) and globally
(considering a relatively long period of time). Research result
evaluation would be reduced to computation, reviewing would become
more time efficient, transparent, some of its current features
such as, for example, anonymity, would become redundant. It may
allow to rigorously compare and uniformize different areas,
projects and activities of mathematics. It would be a quite
helpful research tool for working mathematicians. Correctness and
originality checking may also benefit from a proof
coordinatization since proofs would be written in a more formal
metalanguage. Information about all known mathematical results
could be stored as a single database or similar structure although
this program is not intended to be just another formalization
attempt. Results having low value or low proof complexity should
be marked as such.

A rigorous theory for evaluation of research results must be
developed and implemented regardless of progress of the proposed
more general program since such an evaluation theory would be
valuable in its own right.

\subsubsection{Publishing and social issues - indices and research merits}

The existing system of internationally respected peer reviewed
journals which effectively guides careers of mathematicians seems
to be competition based. Competition is inevitable and productive,
it ensures steady progress and reasonable fairness. Neverheless,
rigorous evaluation of value and complexity of research results,
rigorous publication standards would make the publishing system
more effective and open.

It does not make much sense to use mathematics to precisely
compute various indices such as the h-index, impact factors and
university ratings using well defined formulas if scientific
results are published based upon nontransparent emotional,
business-based or cronyism-based opinions of editors estimating
their values and there are no strict referencing guidelines.
Therefore, in the author's opinion, most ``indices'' which are
supposed to measure scientific productivity are examples of
social-oriented and potentially professionally harmful arguments
serving social needs (such as social ranking struggle, national
needs etc.) more than needs of research and scientific progress.
Formation of individual and collective merits based on social
competition together with deterioration of natural information
processing skills due to computerization may lead to stagnation
and regress.

Mathematical culture should change and the mathematical community
will need to accept the fact that the motivations, logic, the
value of results and the complexity of proofs can be measured.


\subsection{Applications in education}\

Teaching and learning of mathematical (and most other) concepts
and activities are important processes which also need to be
mathematically modelled and analyzed. Following development of
mathematical thinking throughout at least formative years should
be an important part of professional duties of mathematical
community.

\subsubsection{School education}
Possible advances such as a proof complexity measuring, criterions
of simplicity and new encodings used in mathematics could
modernise mathematical studies and make them more effective,
develop optimal learning paths, introduce better or missing
mathematical concepts and activities. Networks of mathematical
concepts and skills taught at school should be reviewed and
updated to correspond also to current advances of computing and
communication technologies which are available to students.

Apart from measuring complexity and weights of proofs research
could also be conducted to measure how easy or difficult
(absolutely or relatively) a definition, an implication, a proof
or other mathematical activity is to understand psychologically or
perform for learners of various abilities and background. Separate
or related learning, teaching and problem solving coordinatization
theories may be necessary to reform education and make it optimal,
truly differentiated and individualized. Below we give two
examples of possible directions of development.

\paragraph{Ideas for dependence based learning and teaching coordinatization theories} For any study course or program we
can define a directed \sl learning graph $\Gamma$:\rm\ the
vertices are knowledge units such as definitions, facts or skills;
a weighted directed edge $a\stackrel{w}{\rightarrow} b$ means that
$b$ must be taught after $a$, the learning difficulty is encoded
in the edge weight $w$. Graph-theoretical considerations given in
\ref{g} can be modified for learning graphs. Given a learning
graph $\Gamma$ we can define for any student $x$ her/his \sl
$\Gamma$-profile (knowledge profile)\rm\ $\pi_{x,\Gamma}$ : a
function from the vertex set $V(\Gamma)$ to a suitable set $K$
which assigns to every knowledge unit $v\in V(\Gamma)$ the level
of knowledge $\pi_{x,\Gamma}(v)$ (an element in $K$) the student
$x$ has with respect to $v$, thus $\pi_{x,\Gamma}\in
\mathcal{F}un(V(\Gamma),K)$, the knowledge profile can also be
interpreted as a set of weights for the $\Gamma$-vertex set - the
weight $\pi_{x,\Gamma}(v)$ of $v$ describes the knowledge or
competence of $v$ which $x$ has. Given a knowledge profile
$\pi_{x,\Gamma}$ we can individually design optimal further
teaching/learning steps. Course goals may be defined in terms of
knowledge profiles, these goals may be defined both individually
and collectively. For any teaching activity $\tau$ we can consider
its impact $f_{\tau}$ on a knowledge profile - $f_{\tau}$ is an
endofunction on $\mathcal{F}un(V(\Gamma),K)$. Similar ideas could
be applied to both natural sciences and other disciplines such as
languages.

\paragraph{Ideas for competence based problem solving theories}
Given a problem $P$ of mathematics or some other discipline we can
define a \sl $P$-solving\rm\ graph $\Gamma$ as follows: the
vertices are states of the solution process, the edges denote
transitions between these states weighted by the skills which are
necessary to perform these transitions. Such graph models would
allow to investigate typical problems used in a given course or
other learning unit, analyze the necessary facts and skills that
need to be taught. They also may be helpful teaching nontrivial
(olympiad) problem solving skills.

\subsubsection{Higher education} University mathematics study programs can also
benefit from a comprehensive proof coordinatization and complexity
theory since such a theory may bring a significant revision,
innovative structuring and modernization of the body of
mathematical knowledge accesible to university students and
graduate students. This revision may affect lists of mathematical
areas, lists of concepts and main facts, course layout designs,
theorems and proof techniques. Proofs and methods which are easy
to understand or proofs containing typical implications could be
collected as a spanning tree for the given theory. The edges
should have uniform complexity. The spanning tree proofs could be
organized into a study course or program. The body of known facts
of every domain should be clearly structured showing the most
complex, "terminal" results.

\subsection{Possible future development and some problems}

\subsubsection{Epistemological problems} A successful progress of
the longterm program discussed in this paper may allow and need
philosophical interpretations.

\paragraph{Methodological Cartesian scepticism and conservativism} A desire to model
and coordinatize mathematical implications, proofs and creativity
contradicts our everyday experience and  represents an instance of
methodological scepticism in mathematics which may be close in
spirit to Cartesian doubt. Everyday mental experiences make most
people believe that implication making is a basic and
unquestionable mental activity, creative mathematical thinking
(and human thinking in general) can not be modelled/mechanised,
its need and modifications can not be discussed. The proposed
research project would challenge these beliefs. Researchers
involved in this project will have to keep this in mind and expect
a great deal of conservativism.

\paragraph{A mathematical justification/regress step} Coordinatization of
mathematical implications and creativity involves mappings of
mathematical proofs to simpler (e.g. discrete, algebraic or
topological) mathematical objects. As it was noted above this
epistemological objective may also be viewed as a generalization
of regress (justification) steps in the standard philosophical
sense, see Pollock (1975). In this program logical implications
which are creative mathematical acts are supposed to be explained
(justified, in philosophical language) in terms of simpler
mathematical objects. Reduction of an applied or pure mathematical
problem or a model to a simpler mathematical object often happens
in mathematics, it has similarities with the standard mathematical
construction of quotient objects.  Since the research in this
program has not even started it may be too early to speculate
about philosophical problems related to the regress step discussed
here such as the mathematical ``problem of the criterion'', see
Chisholm (1989), Cling (2014), or the M\"{u}nchhausen's
(Agrippa's) trilemma, see Albert (1991). The M\"{u}nchhausen
trilemma case determination (i.e. whether the proposed
intra-mathematical regress is cyclic, infinite non-cyclic or
finite) seems to be an important problem.

\subsubsection{Some concrete proposals} We can formulate a few
specific initial research proposals: 1) analyze the body of facts
of Euclidean geometry with respect to the implication modelling
and Hilbertian simplicity idea, create a database of all
nonequivalent logical steps, 2) analyze the body of combinatorics
with respect to structural induction, create a database of all
nonequivalent induction arguments, 3) analyze the body of graph
theory with respect to the proof bundle idea, 4) develop theories
for measuring logical complexity (Hilbertian simplicity) of
standard computational tasks such as solving linear or polynomial
systems of equations, 5) classify invariants and object properties
in a mathematical domain such as, for example, graph theory, with
respect to computational complexity (e.g. polynomial or
NP-complete) of decision problems, study the network of polynomial
reductions, 6) build a map relating theorems and reductions of
computational classes for a given domain, 7) introduce measures of
cog\-ni\-ti\-ve complexity of mathematical activities in school
mathematics courses.

\section{Conclusion} We have given a number of arguments which
justify, encourage and describe a proposal for possible future
research in mathematical logic which can be defined as faithful
mathematical representation of proofs and theories. This can be
called coordinatization of mathematical implications and proofs,
or more generally, mathematical results. The main argument is a
possibility to formalize, map into simpler mathematical objects
and measure mathematical creativity, to separate implication
making from semantic and syntactic issues, to make nontrivial and
creative mathematical theorem proving a computation. Another
argument is a possibility to rigorously measure mathematical
results and to guide the mathematical research in a rigorous and
optimal way. The mathematical culture would greatly benefit from
rigorous standards for mathematical research publications and
other texts. New encoding paradigms for the mathematical language
and mo\-der\-nized mathematical education may also be generated by
modeling of mathematical thinking and introducing coordinates in
the space of proofs - extending the Pythagorean paradigm to
thinking processes.


\end{document}